\documentclass[11pt]{amsart}
\usepackage{amsmath,amsthm,amssymb,amsfonts,amscd}
\usepackage{amsxtra}
\usepackage{eucal}
\usepackage{mathrsfs}
\usepackage{color, graphics}
\usepackage[all]{xy}
\usepackage{hyperref}
\usepackage{pstricks,pst-plot}

\newcommand \Proj {\ensuremath{\mathrm{Proj}}}
\newcommand \im   {\ensuremath{\mathrm{im}}}
\newcommand \Bl {\ensuremath{\mathrm{Bl}}}
\newcommand \mult {\ensuremath{\mathrm{mult}}}
\newcommand \ini {\ensuremath{\mathrm{in}}}
\newcommand \Tor {\ensuremath{\mathrm{Tor}}}
\newcommand \Trisec {\ensuremath{\mathrm{Trisec}}}

\newcommand \Sec {\ensuremath{\mathrm{Sec}}}
\newcommand \Tan {\ensuremath{\mathrm{Tan}}}
\newcommand \Sym {\ensuremath{\mathrm{Sym}}}

\newcommand \codim {\ensuremath{\mathrm{codim}}}
\newcommand \coker {\ensuremath{\mathrm{coker}}}
\newcommand \length{\ensuremath{\mathrm{length}}}

\newcommand \depth {\ensuremath{\mathrm{depth}}}

\def\P{{\mathbb P}}
\def\Z{{\mathbb Z}}

\newcommand \krm[1] {K_{#1}^R(M)}
\newcommand \cone[1] {C_{#1}(\overline{\varphi})}
\newcommand \kosv[3] {\wedge^{#1} V\otimes {#2}_{#3}}
\newcommand \kosw[3] {\wedge^{#1} W\otimes {#2}_{#3}}

\newcommand \st[1] {\stackrel{#1}{\longrightarrow}}

\newcommand \reg {\mathrm {reg}}

\newcommand \lra {\rightarrow}
\newcommand \llra {\longrightarrow}
\newcommand \p {\partial}

\def\ds{\displaystyle}

\theoremstyle{plain} 
\newtheorem{Thm}{Theorem}[section]
\newtheorem{Prop}[Thm]{Proposition}
\newtheorem{Coro}[Thm]{Corollary}
\newtheorem{Lem}[Thm]{Lemma}

\theoremstyle{definition}
\newtheorem{Def}[Thm]{Definition}
\newtheorem{Ex}[Thm]{Example}
\newtheorem{Qu}[Thm]{Question}
\newtheorem{Remk}[Thm]{Remark}

\newtheorem{Clm}[Thm]{Claim}

\newtheorem{Nota}[Thm]{Notation}
\numberwithin{equation}{section}

\begin{document}

\title{Graded mapping cone theorem, multisecants and syzygies}
\author[J.\ Ahn and S.\ Kwak]{Jeaman Ahn and Sijong Kwak}
\address{Department of Mathematics, Korea Advanced Institute of Science and Technology,
373-1 Gusung-dong, Yusung-Gu, Daejeon, Korea}
\email{jeamanahn@kongju.ac.kr}
\address{Department of Mathematics, Korea Advanced Institute of Science and Technology,
373-1 Gusung-dong, Yusung-Gu, Daejeon, Korea}
\email{skwak@kaist.ac.kr}
\thanks{The second author was supported in part by Basic Science Research Program through 
the National Research Foundation of Korea(NRF) funded by the Ministry of Education, Science and Technology
(grant no.2009-0063180) and KRF(grant No.2005-070-C00005).}

\begin{abstract}
Let $X$ be a reduced closed subscheme in $\mathbb P^n$.
As a slight generalization of property $\textbf{N}_p$ due to Green-Lazarsfeld,
we can say that $X$ satisfies property $\textbf{N}_{2,p}$
scheme-theoretically if there is an ideal $I$ generating the ideal sheaf $\mathcal I_{X/\P^n}$
such that $I$ is generated by quadrics and there are only linear syzygies up to $p$-th step
(cf. \cite{EGHP1}, \cite{EGHP2}, \cite{V} ). Recently, many algebraic and geometric results have been proved for
projective varieties satisfying property $\textbf{N}_{2,p}$(cf. \cite{CKP}, \cite{EGHP1}, \cite{EGHP2} \cite {KP}).
In this case, the Castelnuovo regularity and normality can be obtained by the blowing-up method
as $\reg(X)\le e+1$ where $e$ is the codimension of a smooth variety $X$ (cf. \cite{BEL}).
On the other hand, projection methods have been very useful and powerful in bounding Castelnuovo
regularity, normality and other classical invariants in geometry(cf. \cite{BE}, \cite{K}, \cite{KP},
\cite{L} \cite {R}).

In this paper, we first prove the graded mapping cone theorem on partial eliminations as a general algebraic tools
and give some applications.
Then, we bound the length of zero dimensional intersection of $X$ and a linear space $L$
in terms of graded Betti numbers and deduce a relation between $X$ and its projections
with respect to the geometry and syzygies in the case of projective schemes satisfying property $\textbf{N}_{2,p}$
scheme-theoretically.
In addition, we give not only interesting information on the regularity of fibers and multiple loci
for the case of $\textbf{N}_{d,p}, ~~d\ge 2$ but also geometric structures for projections according to moving the center.
\\
\end{abstract}

\maketitle
\noindent 2000 Mathematics Subject Classification: 14N05, 13D02, 14M17.\\
Keywords: linear syzygies, graded mapping cone, Castelnuovo-Mumford regulaity, partial elimination ideal.

\tableofcontents \setcounter{page}{1}

\section{Introduction}\label{section_1}
Let $X$ be a non-degenerate reduced closed subscheme in $\P(V)$ with the saturated ideal
$I_X=\bigoplus_{m=0}^{\infty}H^0(\mathcal I_{X/\P^n}(m))$ where
$V$ is an $(n+1)$-dimensional vector space over an algebraically closed field $k$ of characteristic zero
and $R=k[x_0,\ldots,x_n]$ be the coordinate ring of $\P(V)$.
As Eisenbud et al in \cite{EGHP1} introduced the notion
$\textbf {N}_{d,p}$ for some $d\geq 2$, we say that $X$(or $I_X$) satisfies the property $\textbf {N}_{d,p}$ if $\Tor^R_i(R/I_X,k)$
is concentrated in degree $<d+i$ for all $i\leq p$, which is equivalent to the condition that the truncated ideal
$(I_X)_{\geq d}$ is generated by homogeneous forms of degree $d$ and has a linear resolution until the first $p$
steps.

The case of $d=2$ has been of particular interest and there have been many classical conjectures and known
results for highly positive embeddings and the canonical embedding of a smooth variety $X$.
The property $\textbf{N}_{2,p}$ is the same as the property $\textbf {N}_p$ defined by
Green and Lazarsfeld if $X$ is projectively normal. The property $\textbf {N}_{2,1}$ means that $I_X$ is generated
by quadrics and the property $\textbf {N}_{2,2}$ means that there are only linear relations on quadrics
in addition to the property $\textbf {N}_{2,1}$.
In \cite{EGHP1} the authors have exhibited some cases in which there is an interesting connection
between the minimal free resolution of $I_X$ and the minimal free resolution of $I_{X\cap L, L}$
where $L$ is a linear subspace of $\P(V)$.
They have shown that if $X$ satisfies $\textbf{N}_{2,p}$ and $\dim X\cap L \leq 1$ for some linear space
$L$ of dimension $\leq p$, then $I_{X\cap L, L}$ is $2$-regular. They also gave the conditions
when the syzygies of $X$ restrict to the syzygies of the intersection.\\

As a slight generalization, we can say that $X\subset \mathbb P(V)$ satisfies property $\textbf{N}_{2,p}$
scheme-theoretically if there is an ideal $I$ generating the ideal sheaf $\mathcal I_{X/\P^n}$
such that $I$ is generated by quadrics and there are only linear syzygies up to $p$-th step
(cf.~\cite{EGHP1}, \cite{V}). For example, if $I_X$ satisfies property $\textbf{N}_{2,p}$ then a
general hyperplane section $X\cap H$ satisfies property $\textbf{N}_{2,p}$ scheme-theoretically
because $\frac{I_X+(H)}{(H)}$ is generated by quadrics in $R/(H)$ and has only linear syzygies up to first $p$-th steps
even if $\frac{I_X+(H)}{(H)}$ is not saturated in general.
If $X$ is smooth and is cut out by quadrics scheme-theoretically, then the Castelnuovo regularity and
normality are easily obtained by the blowing-up method and Kawamata-Viehweg vanishing theorem. In particular,
$\reg(X)\le e+1$ where $e$ is the codimension of $X$ (cf. \cite{BEL}).

On the other hand, projection methods have been very useful and powerful in bounding the Castelnuovo
regularity, higher order normality and other classical invariants in geometry(cf.~\cite{BE}, \cite{K}, \cite{KP},
\cite{L}, \cite {R}).
Consider a linear projection $\pi_{\Lambda}: X\to Y\subset \P^{n-t}=\P(W), ~W\subset V$ from the center $\Lambda=\P^{t-1}$
such that $\Lambda\cap X=\emptyset$. What can we say about a connection between the minimal free resolution of $I_X$
and the minimal free resolution of $I_Y$? We are mainly interested in homological, cohomological,
geometric and local properties of projections as the center moves in an ambient space.
Note that for graded ideals or modules which are not saturated, the Koszul techniques of Green and Lazarsfeld are
not so much adaptable to understand their syzygies.

Our basic idea is to compare the graded Betti numbers of $R/I_X$ as a $R$-module and those of $R/I_X$ as a $\Sym(W)$-module
via the graded mapping cone theorem and to interpret its geometric meanings.

The paper is organized as follows: in Section 2, we prove the graded mapping cone theorem on partial eliminations
and give some applications. When applying this theorem to the projection $\pi_{\Lambda}: X\twoheadrightarrow Y\subset \P^{n-t}$
from the center $\Lambda=\P^{t-1}$ such that $\Lambda\cap X=\emptyset$, we prove that every fiber of {\it arbitrary}
projection $\pi_{\Lambda}:X\to Y$ is $(d-1)$-normal if $X$ satisfies the property
$\textbf {N}_{d,p}$ scheme-theoretically and $\dim\Lambda < p$ and recover some results in \cite{EGHP1}:
a linear section $X\cap L$ is $d$-regular if $\dim(X\cap L)=0$ and $\dim L\le p$. In particular, a projective variety satisfying property
$\textbf{N}_{2,p}$ scheme-theoretically has no $(p+2)$-secant $p$-plane. As a generalization, we bound the possible maximal length
of $X\cap L$ in terms of the graded Betti numbers (see Theorem~\ref{Theo:1-3}).

In Section 3, we study the effects of the property $\textbf{N}_{2,p}$ scheme-theoretically on the Castelnuovo normality
and defining equations of projected varieties (see Theorem ~\ref{Main result N_{2,p}} and Theorem ~\ref{birational N_{2,p}}).
Using the partial elimination ideal theory due to M. Green \cite{G},
we give some information on the multiple loci for birational projections (Theorem ~\ref{N_{d,p}-geometry}) for the case of the  property
$\textbf{N}_{d,2}$. Moreover, we obtain the relation between the regularity of the $i$-th partial elimination ideals for all $i\ge 1$ and
syzygies of projections even though the Castelnuovo normality is very delicate and difficult to control under projections.

In Section 4, we deal with some properties of projections, e.g. the number of quadratic equations and the depth of projected varieties
according to moving the center~~(e.g., Proposition ~\ref{prop_moving}).
In particular, for birational projections, we show that the singular locus of the projected variety is
a linear subspace for $p\ge 2$. From this fact, we give some interesting non-normal examples and applications. \\

{\bf Acknowledgements} It is our pleasure to thank the referee for valuable comments and further directions.
The second author would like to thank Korea Institute of Advanced Study(KIAS) for supports and hospitality during his stay there
for a sabbatical year.
\section{Graded mapping cone theorem and applications}\label{section_2}
The mapping cone under projection and its related long exact sequence is our starting point to understand
algebraic and geometric structures of projections.\\
\begin{itemize}
 \item $W=\bigoplus^{n}_{i=1}k\cdot x_i\subset V=\bigoplus^{n}_{i=0}k\cdot x_i$ : vector spaces over k.
 \item $S_1=k[x_1,\ldots,x_n]\subset R=k[x_0,\ldots,x_n]$ : polynomial rings.
 \item $M$ : a graded $R$-module (which is also a graded $S_1$-module).
 \item $K_{\bullet}^{S_1}(M)$ : the graded Koszul complex of $M$ as follows:
 $$0\rightarrow \wedge^n W\otimes {M}\rightarrow\cdots\rightarrow\wedge^2 W\otimes {M}\rightarrow W\otimes {M}\rightarrow M\rightarrow 0$$ whose graded components are $K_{i}^{S_1}(M)_{i+j} = \wedge^i W\otimes {M_j}$.
\end{itemize}
Consider the multiplicative map $\varphi: M(-1) \stackrel{\times x_0}{\longrightarrow} M$ as a graded $S_1$-module homomorphism
such that $\varphi(m)=x_0\cdot m$. Then we have the induced map
$$\overline{\varphi}: \Bbb F_{\bullet}=K_{\bullet}^{S_1}(M(-1)) \stackrel{\times x_0}{\longrightarrow} \Bbb G_{\bullet}=K_{\bullet}^{S_1}(M)$$
between graded complexes.

Now, we construct the mapping cone $(C_{\bullet}(\overline{\varphi}),\p_{\,\,\overline{\varphi}})$ induced by $\overline{\varphi}$
such that $C_{\bullet}(\overline{\varphi})=\Bbb G_{\bullet}\bigoplus\Bbb F_{\bullet}[-1]$ and
\begin{itemize}
\item $C_i(\overline{\varphi})_{i+j}=[\Bbb G_{i}]_{i+j}\bigoplus
[\Bbb F_{i-1}]_{i+j}=\wedge^i W\otimes M_j\oplus \wedge^{i-1} W\otimes M_j$.\\
\item the differential $\p_{\,\,\overline{\varphi}}: C_i(\,\,\overline{\varphi})
\rightarrow C_{i-1}(\overline{\varphi})$ is given by
\[\p_{\,\,\overline{\varphi}}=\left(%
\begin{array}{cc}
  \p & \overline{\varphi} \\
  0 & -\p \\
\end{array}%
\right),\] where $\p$ is the differential of Koszul complex $K_{\bullet}^{S_1}(M)$.
\end{itemize}
Finally, the mapping cone $(C_{\bullet}(\overline{\varphi}),\p_{\,\,\overline{\varphi}})$
becomes a complex over $S_1$ and we have the exact sequence of
complexes
\begin{equation}\label{eq:101}
0\longrightarrow \Bbb G_{\bullet} \longrightarrow C_{\bullet}(\overline{\varphi}) \longrightarrow \Bbb F_{\bullet}[-1]\longrightarrow 0.
\end{equation}

From the exact sequence (\ref{eq:101}), we have a long exact
sequence in homology:\\
\begin{equation}\label{eq:102}
\begin{array}{ccccccccccccc}
\st{} \Tor_{i}^{S_1}(M,k)_{i+j}& \st{}
H_i(C_{\bullet}(\overline{\varphi}))_{i+j}&\st{}&&
\\[2ex]
&\Tor_{i-1}^{S_1}(M,k)_{i+j-1}&\st{\delta = \times x_0}&
\Tor_{i-1}^{S_1}(M,k)_{i+j}& \st{} &
\end{array}
\end{equation}
and the connecting homomorphism $\delta$ is the multiplicative map induced by
$\overline{\varphi}$.

In the following Lemma ~\ref{lem:101}, we claim that $\Tor^R(M,k)$ can be obtained by the homology of the mapping cone.
\begin{Lem}\label{lem:101}
Let $M$ be a graded $R$-module. Then we have the following natural
isomorphism:
$$\Tor_i^R(M,k)_{i+j}\simeq H_i(\cone{\bullet})_{i+j}.$$

\end{Lem}
\begin{proof}

Let $\krm{\bullet}$ be the Koszul complex of a graded $R$-module
$M$. Then the graded component in degree $i+j$ of $\krm{i}$ is $\krm{i}_{i+j}=\kosv{i}{M}{j}$.
Note that $\wedge^{i}V\cong
[x_0\wedge(\wedge^{i-1}W)]\oplus\wedge^{i}W$. Hence we see that
the Koszul complex $\krm{i}$ has the following canonical
decomposition in each graded component:\\
\begin{equation}\label{eq:103}
\begin{array}{cccccccccccccccccccc}
&&\ds\kosw{i}{M}{j}&&&\\[1ex]
\ds \krm{i}_{i+j}&\ds \cong& \ds\bigoplus&\ds\cong&\cone{i}_{i+j}.\\[1ex]
&&[x_0\wedge(\wedge^{i-1} W)]\otimes {M}_{j}&\\[1ex]
\end{array}
\end{equation}
Using the decomposition (\ref{eq:103}), we can verify that the
following diagram is commutative:\\
\begin{equation}\label{eq:104}
\begin{array}{cccccccccccccccccccc}
\ds \krm{i}_{i+j}&\ds \st{\ds \cong} & \cone{i}_{i+j}\\[2ex]
\Big\downarrow\vcenter{\rlap{$\p$}}&&\Big\downarrow\vcenter{\rlap{$\p_{\overline{\varphi}}$}} \\[2ex]
\ds \krm{i-1}_{i+j}&\ds \st{\ds \cong} & \cone{i-1}_{i+j}.\\[2ex]
\end{array}
\end{equation}
Therefore, we have a natural isomorphism
$\Tor_i^R(M,k)_{i+j}\simeq H_i(\cone{\bullet})_{i+j}$.

\end{proof}

From the long exact sequence (\ref{eq:102}) and
Lemma~\ref{lem:101}, we obtain the following useful Theorem.

\begin{Thm}\label{Theo:1-2} Let $S_1=k[x_1,\ldots,x_n]\subset R=k[x_0,x_1\ldots,x_n]$ be
polynomial rings. For a graded $R$-module $M$, we have the
following long exact sequence:
\begin{equation*}
\begin{array}{ccccccccccccc}
\llra\Tor_{i}^{S_1}(M,k)_{i+j}\llra\Tor_{i}^R(M,k)_{i+j}\llra
\Tor_{i-1}^{S_1}(M,k)_{i+j-1}\llra&  &&
\\[2ex]
\st{\delta}\Tor_{i-1}^{S_1}(M,k)_{i+j}\llra\Tor_{i-1}^R(M,k)_{i+j}\llra
\Tor_{i-2}^{S_1}(M,k)_{i+j-1}\st{\delta}&\cdots
\end{array}
\end{equation*}
whose connecting homomorphism $\delta$ is the multiplicative map $\times\, x_0$.
\end{Thm}\begin{proof}
It is clear from (\ref{eq:102}) and
Lemma~\ref{lem:101}.
\end{proof}

Note that Theorem~\ref{Theo:1-2} gives us an useful information about
syzygies of outer projections (i.e. isomorphic or birational projections)
of projective varieties.

As a first step, we obtain the following important Corollary.

\begin{Coro}\label{Coro:101}
Let $I\subset R$ be a homogeneous ideal such that $R/I$ is a
finitely generated $S_1$-module. Assume that $I$ admits $d$-linear
resolution as a $R$-module up to $p$-th step for $p\geq 2$. Then, for $1\leq i\leq
p-1$,
\begin{enumerate}
\item[(a)] the minimal free resolution of $R/I$ as a graded
$S_1$-module is
$$\rightarrow L_{p-1}\rightarrow\cdots\rightarrow L_1
\rightarrow S_1(-d+1)\oplus \cdots\oplus S_1(-1)\oplus S_1\rightarrow R/I \rightarrow
0,$$
where $L_i=S_1(-d+1-i)^{\oplus \beta^{S_1}_{i,d-1}}$ for all $1\le i\le p-1$;
\item[(b)] in particular, \,\,$\beta^{S_1}_{i,d-1}=(-1)^{i}+\sum_{1\leq j\leq
 i}(-1)^{j+i}\beta_{j,d-1}^R(R/I)$.
\end{enumerate}
\end{Coro}
\begin{proof}
(a)
First, consider the exact sequence
\begin{equation*}
\begin{array}{ccccccccccccc}
\lra \Tor_{1}^R(R/I,k)_{j}& \lra
\Tor_{0}^{S_1}(R/I,k)_{j-1}&\st{\delta}&&
\\[2ex]
&\Tor_{0}^{S_1}(R/I,k)_{j}&\lra& \Tor_{0}^R(R/I,k)_{j}&\lra &0.
\end{array}
\end{equation*}
Since $\Tor_{1}^R(R/I,k)_{j}=0$ for all $j\neq d$ and $\Tor_{0}^R(R/I,k)_{j}=0$ for
all $j\neq 0$, we obtain that $\beta^R_{0,0}=\beta^{S_1}_{0,j}=1$ for all
$0\leq j\leq d-1$ and $\beta^{S_1}_{0,j}=0$ for all $j\notin
\{0,1,\ldots,d-1\}$ from the finiteness of $R/I$ as a $S_1$-module.

Note that $\Tor_{i}^R(R/I,k)_{i+j}=0$ for $1\le
i\le p$ and $j\neq d-1$ by assumption that $I$ is $d$-linear up to
$p$-th step. Applying Theorem~\ref{Theo:1-2} for $M=R/I$, we have
an isomorphism induced by $\delta= \times\, x_0$
\[\Tor_{i-1}^{S_1}(R/I,k)_{(i-1)+j}\st{\delta}\Tor_{i-1}^{S_1}(R/I,k)_
{(i-1)+(j+1)}\] for $1\le i\le p$ and for all $j\notin \{d-2,d-1\}$.
Hence we conclude that
$$\Tor_{i-1}^{S_1}(R/I,k)_{(i-1)+j}=0 \,\,\text { for } 2\le i\le p \,\,\text {  and  }
\,\,j\neq d-1  $$ since $R/I$ is finitely generated as an $S_1$-module,
which means that
$$L_i=S_1(-d-i+1)^{\oplus \beta^{S_1}_{i,d-1}}\,\,\text { for }
\,\,\,1\le i \le p-1.$$
(b) Note that we have
\begin{equation*}
0 \lra \Tor^{S_1}_{i}(R/I,k)_{i+d-1}\lra
\Tor^R_{i}(R/I,k)_{i+d-1}\lra \Tor^{S_1}_{i-1}(R/I,k)_{i+d-2}\lra
0
\end{equation*}
for $1\leq i\leq p-1$ such that
\[\beta_{i,d-1}^{S_1}(R/I)=\beta_{i,d-1}^R(R/I)-\beta_{i-1,d-1}^{S_1}(R/I).\]
Then, by induction on $p$, we get the desired result.
\end{proof}

\begin{Nota}
In this paper, we use the following notations:
\begin{itemize}
 \item[$\bullet$] $R=k[x_0,\ldots,x_n]=\Sym(V)$ and $S_t=k[x_t,x_{t+1}\ldots,x_n]=\Sym(W)$ are two polynomial rings where
 $W\subset V$, $\codim(W,V)=t$.
 \item[$\bullet$] $\Lambda=\P(U)=Z(x_t,x_{t+1},\cdots, x_n)$ is a linear space in $\P^n$ where $U$ is a $t$-dimensional
     vector space with a basis $\{x_0, x_1, \cdots, x_{t-1}\}$.
 \item[$\bullet$] $X$ is assumed to be a non-degenerate reduced projective scheme if unspecified.
 \item [$\bullet$]$\pi_{\Lambda}$ : $X\to Y_t= \pi_{\Lambda}(X)\subset\mathbb P^{n-t}=\P(W)$ is the projection
     from the center $\Lambda$ and $\Lambda \cap X=\phi$.
 \item[$\bullet$] $\beta^R_{i,j}(M):=\dim_k\Tor_i^{R}(M,k)_{i+j}$ for a finitely generated $R$-module $M$.
 \item[$\bullet$] $H^i_{*}(\mathcal F):=\bigoplus_{\ell \in \Z} H^i(\mathcal F(\ell))$ and $h^i(\mathcal F)=\dim H^i(\mathcal F)$
    for a coherent sheaf $\mathcal F$.
\end{itemize}
\end{Nota}

From now on, we consider a projection $\pi_{\Lambda}$ : $X\to Y_t= \pi_{\Lambda}(X)\subset\P(W)$ where $\dim \Lambda =t-1\ge 0,
\Lambda \cap X=\phi$. Then, the following basic sequence
$$ 0 \longrightarrow R/I_X \longrightarrow E\longrightarrow H^1_{*}(\mathcal{I}_X)\longrightarrow 0 \,\,\,\,\,\,
(\text {as $S_t$-modules} )
$$ is also exact as finitely generated $S_t$-modules as Lemma ~\ref{fin} shows.
Furthermore, it would be very useful to compare their graded Betti tables by the graded mapping cone theorem.

\begin{Lem}\label{fin}
Let $I$ be a homogeneous ideal defining $X$ scheme-theoretically in $\P^n$. Then $R/I$ and $E=\bigoplus_{\ell \in \Z} H^0 (X,\mathcal{O}_X(\ell))$
are finitely generated $S_t$-modules.
\end{Lem}
\begin{proof}
For $0\le i\le t-1$, let $p_i=[0,\ldots,0,1,0,\ldots,0]\in \Lambda$
be the point whose $(i+1)$-th coordinate is 1. Since $p_i\notin X$, for some $m_i>0$, there is a homogeneous polynomial $f_i$ in
$I$ which is of the form $f_i = x_i^{m_i}+ g_i$ where $g_i\in k[x_{0},x_1,\ldots,x_n]$ is a homogeneous polynomial of degree
$m_i$ with the power of $x_i$ less than $m_i$.
Hence $R/I$ is generated by monomials of the form $x_0^{n_0}\cdots x_{t-1}^{n_{t-1}}$, $n_i< m_i$ for all $0\leq i\leq t-1$
as a $S_t$-module. Next, from the exact sequence $ 0 \longrightarrow R/I_X \longrightarrow E\longrightarrow
H^1_{*}(\mathcal{I}_X)\longrightarrow 0 $ as $S_t$-modules,
$E$ is also a finitely generated $S_t$-module.
\end{proof}

\begin{Remk}
For an inner projection of $X$ from the center  $q\in X$, let $Y_1=\overline{\pi_{q}(X)}$
be the Zariski-closure of $\pi_{q}(X)$ in $\mathbb P^{n-1}$. Then $R/I_X$ is {\it not} a finitely generated $S_1$-module.
\end{Remk}

The following theorem is a generalization of Corollary~\ref{Coro:101}, which is related to the existence of multisecant plane
(cf. Theorem~\ref{Theo:1-3}).

\begin{Thm}\label{Prop:201}
Suppose that $X$ satisfies property $\textbf {N}_{d,p}$ scheme-theoretically with an ideal $I$. Consider the linear projection
$\pi_{\Lambda}$:$X\to \P(W), \Sym(W)=S_t$ from the center $\Lambda$ such that $\Lambda \cap X=\phi, \,\,\Lambda=\P(U)=\P^{t-1}$.
Then, $\overline R=R/I_{\ge d}$ has the simplest syzygies up to $(p-t)$-th
step as $S_t$-module for $1\le t\le p$,
\begin{equation}\label{eq:211}
\rightarrow L_{p-t}\rightarrow\cdots \rightarrow L_1 \rightarrow \bigoplus^{d-1}_{i=0} \Sym^i(U)\otimes
S_t(-i)\rightarrow \overline R
\rightarrow  0.
\end{equation}
where $L_i=S_t(-i-d+1)^{\beta_{i,d-1}^{S_t}}$ for $1\le i\le p-t$ and $\Sym^i(U)=H^0(\mathcal O_{\Lambda}(i))$ is
a vector space of homogeneous forms of degree $i$ generated by $U$.

In particular, if $d=2$ then the minimal free resolutions of $R/I$ is
$$\rightarrow S_t(-p+t-1)^{\beta_{p-t,1}^{S_t}}\rightarrow \cdots \rightarrow S_t(-2)^{\beta_{1,1}^{S_t}}\rightarrow S_t
\oplus S_t(-1)^{t}\rightarrow R/I\rightarrow  0.
$$
\end{Thm}
\begin{proof}
Let $S_t=k[x_t,\ldots,x_n]$ be a polynomial ring for $0\leq t\leq p$ and let
\[ \rightarrow L_{p-t}\rightarrow\cdots \rightarrow L_1 \rightarrow L_0\rightarrow \overline R
\rightarrow  0
\]
be the minimal free resolution of $\overline R$ as a $S_t$-module.
We will give a proof by induction on $t\geq 1$. For $t=1$, the result (\ref{eq:211}) follows directly from Corollary~\ref{Coro:101}.
For $t>1$,  by induction hypothesis, we can assume that for  $1\leq \alpha\leq p-(t-1)$
\[\Tor_{\alpha}^{S_{t-1}}(\overline R,k)_{j}=0 \quad \text{ if }j\neq \alpha+d-1.\]
Using an exact sequece by mapping cone theorem for $\alpha\geq 1$
\begin{equation*}\label{eq:212}
\begin{array}{ccccccccccccc}
\llra \Tor_{\alpha+1}^{S_{t-1}}(\overline R,k)_{j}& \llra
\Tor_{\alpha}^{S_{t}}(\overline R,k)_{j-1}&\st{\delta}&&
\\[2ex]
&\Tor_{\alpha}^{S_{t}}(\overline R,k)_{j}&\llra& \Tor_{\alpha}^{S_{t-1}}(\overline R,k)_{j}&\llra &
\end{array}
\end{equation*}
we can also show by similar argument used in Corollary~\ref{Coro:101} that
\[\Tor_{\alpha}^{S_{t}}(\overline R,k)_{j}=0 \quad \text{ if }j\neq \alpha+d-1,\]
equivalently, $L_\alpha=S_t(-\alpha-d+1)^{\beta_{\alpha,d-1}^{S_t}}$ for $1\le \alpha\le p-t$.

It remains to show that $L_0=\bigoplus^{d-1}_{i=0} \Sym^i(U)\otimes
S_t(-i)$. Note that the set $\{\Sym^i(U)|0\le i\le d-1\}$ should be contained in any generating set of
$\overline R$ as a $S_t$-module because there is no relation of degree $\le d-1$ in $\overline R$. So, we have to show that
$\{\Sym^i(U)|0\le i\le d-1\}$ is actually the set of all
generators. This can be done by the dimension counting. Let us
prove this by induction on $t$. In the case of $t=1$, the result
easily follows from Corollary~\ref{Coro:101} (a). If we have $t>1$ then, by
induction hypothesis, we see that for all $i\le d-1$,
$$\dim_k\Tor^{S_{t-1}}_0(\overline R,k)_{i}=\binom{i+t-2}{t-2}, \quad \Tor^{S_{t-1}}_1(\overline R,k)_{i}=0$$
and we have the following sequence from the mapping cone construction
\[0\rightarrow \Tor^{S_t}_0(\overline R,k)_{i-1}\rightarrow \Tor^{S_t}_0(\overline R,k)_i\rightarrow
\Tor_{0}^{S_{t-1}}(\overline R,k)_i\rightarrow 0,\] for each $0\leq i<d$. Hence, we obtain
\[
\dim_k\Tor_{0}^{S_{t}}(\overline R,k)_i=\sum_{m=0}^{i}\binom{m+t-2}{t-2}=\binom{i+t-1}{t-1},
\]
as we wished.
\end{proof}

\begin{Def}[\cite{KP}]
Consider three vector spaces $W\subset V\subset H^0(X,\mathcal{O}_X(1))$ and suppose that $t=\mbox{codim}(W,V)$ and
$\alpha=\codim(W,H^0(\mathcal{O}_X(1)))$. We say that $R/I_X$ (resp. $E$) satisfies property $\textbf {N}^{S}_p$ if $R/I_X$ (resp. $E$)
have the simplest minimal free resolutions until $p$-th step as graded $S_t$-modules;
\begin{equation}
\cdots \rightarrow E_p \rightarrow E_{p-1}\rightarrow \cdots\rightarrow E_1 \rightarrow S_t\oplus S_t(-1)^{\oplus \alpha}
\rightarrow E \rightarrow 0
\end{equation}
where $E_i= S_t(-i-1)^{\oplus \beta_{i,1}}$ for $1\le i\le p$ \,\, and
\begin{equation}
\cdots\rightarrow \tilde{L}_p \rightarrow\tilde{L}_{p-1}\rightarrow \cdots \rightarrow \tilde{L}_1
\rightarrow S_t\oplus S_t(-1)^{\oplus \,\,t}\rightarrow R/I_X\rightarrow 0
\end{equation}
where $\tilde{L}_i=S_t(-i-1)^{\oplus \tilde\beta_{i,1}}, 1\le i\le p$.
\end{Def}

On the other hand, we have the similar result for $E=\bigoplus_{\ell \in \Z} H^0 (X,\mathcal{O}_X (\ell))$ as the following
proposition shows.
\begin{Prop}\label{Prop:202}
In the same situation as in the Theorem ~\ref{Prop:201}, suppose $E$(or $R/I_X$) satisfies property $\textbf {N}^S_{p}$ as $R$-module for some
$p\geq 2$. Then $E$ (or $R/I_X$) also satisfies property $\textbf {N}^S_{p-t}$ as $S_t$-module under the projection morphism $\pi_{\Lambda}$ : $X\to Y_t\subset
\mathbb P^{n-t}=\P(W), 1\le t \le p.$ \end{Prop}
\begin{proof} When $t=1$, we can similarly show that $E$ satisfies property
$\textbf {N}^S_{p-1}$ as an $S_1$ module by using Theorem~\ref{Theo:1-2} for $M=E$ and the vanishing $\beta_{i,j}^R(E)=0, 0\le i\le p, j\ge 2$.
As a consequence, $E$ has the following simplest resolution;
\begin{equation}
\cdots \rightarrow E_{p-1} \rightarrow E_{p-2}\rightarrow \cdots
\rightarrow E_1 \rightarrow S_1\oplus S_1(-1)^{\oplus
\alpha}\rightarrow E \rightarrow 0
\end{equation}
where $E_i=S_1(-i-1)^{\oplus \beta_{i,1}}$ for $1\le i\le p-1$. For $t\ge 2$, letting $S_i=k[x_i,x_{i+1}\ldots,x_n]$, we inductively check
that if $E$ satisfies property $\textbf {N}^S_{p-i}$ as an $S_{i}$-module, then $E$ also satisfies property $\textbf {N}^S_{p-i-1}$ as an $S_{i+1}$-module
by the same argument as in the Theorem ~\ref{Prop:201}. For $R/I_X$, the proof is exactly same.
\end{proof}

\begin{Remk}
For isomorphic projections, the above Proposition~\ref{Prop:202} is in fact a simple algebraic reproof of Theorem 2 in \cite {CKP},
Theorem 1.2 in \cite{KP}, and for birational projections, see a part of Theorem 3.1 in \cite{P2}. Indeed, for any regular projection
$\pi_{\Lambda}$ : $X\to Y_t\subset \mathbb P^{n-t}=\P(W), 1\le t \le p$, there is an exact sequence:
$
0\rightarrow \Tor_i^{S_t}(E,k)_{i+j}\rightarrow H^1({\wedge}^{i+1}\mathcal M\otimes {\pi_{\Lambda}}_{*}\mathcal O_X(j-1))
\st{\alpha_{i,j}}{\wedge}^{i+1}W\otimes H^1({\pi_{\Lambda}}_{*}\mathcal O_X(j-1))\rightarrow \cdots
$. From the following commutative diagram:
\begin{equation*}\label{eq:104}
\begin{array}{cccccccccccccccccccc}
\ds H^1(Y_t, {\wedge}^{i+1}\mathcal M_W\otimes {\pi_{\Lambda}}_{*}\mathcal O_X(j-1))& \st{\ds \alpha_{i,j}} & {\wedge}^{i+1}W\otimes
H^1({\pi_{\Lambda}}_{*}\mathcal O_X(j-1))\\[2ex]
\parallel &&\parallel\\[2ex]
H^1(X, {\wedge}^{i+1}{\pi}_{\Lambda}^*\mathcal M_W\otimes\mathcal O_X(j-1))& \st{\ds \widetilde \alpha_{i,j}} &{\wedge}^{i+1}W\otimes
H^1(\mathcal O_X(j-1)),\\[2ex]
\end{array}
\end{equation*}
it was shown (cf. \cite {CKP},\cite{KP}, \cite{P2}) that $\widetilde \alpha_{i,j}$ is injective by induction for all $1\le i\le p-t$
and $j\ge 2$. Thus, $E$ satisfies property $N^S_{p-t}$ as $S_t$-module.  \qed \end{Remk}

The following theorem gives us a geometric meaning of property $N_{d,p}$ with respect to multisecant planes. Note that part (b) was
also proved in Theorem 1.1 in \cite {EGHP1} with a different method.

\begin{Thm}\label{Theo:1-3}
Suppose that $X$ satisfy property $\textbf {N}_{d,p}$ scheme-theoretically in $\P^n$. Consider
the projection $\pi_{\Lambda}$:$X\to Y_t\subset \mathbb P^{n-t}=\P(W)$ from the center
$\Lambda$ such that $\Lambda \cap X=\phi\,, \,\,\Lambda=\P(U)=\P^{t-1}, t\le p$.
Then,
\begin{itemize}
\item[(a)]
every fiber of $\pi_{\Lambda}$ is $(d-1)$-normal, i.e. $\reg(\pi_{\Lambda}^{-1}(y))\le d$ for all $y\in Y_t$.
So, $\length(\pi_{\Lambda}^{-1}(y))\le \binom{t+d-1}{t}$;
\item[(b)] $\reg(X\cap L)\le d$ for any linear section $X\cap L$ as a finite scheme where $L=\P^{k_0}$, $1\le k_0\le p$.
In particular, for a projective variety satisfying property $\textbf {N}_{2,p}$, there is no $(p+2)$-secant $p$-plane.
\item[(c)] Suppose that $X$ satisfies $\textbf {N}_{2,p}$ and $\textbf {N}_{3,p+1}$ scheme-theoretically for $p\ge 0$.
If there is a $l$-secant $(p+1)$-plane then $$l\leq p+2+\min\{p+1, \beta^R_{p+1,2}(R/I)\}.$$
\end{itemize}
\end{Thm}

\begin{proof} Choose an ideal $I$ defining $X$ with $\textbf {N}_{d,p}$ scheme-theoretically.
For a proof of (a), consider the minimal free resolution of $R/(I)_{\ge d}$ in Proposition ~\ref{Prop:201}. (b), namely,
$$\cdots\rightarrow S_t(-d)^{\beta_{1,d-1}^{S_t}}\rightarrow \bigoplus^{d-1}_{i=0}\Sym^i(U)\otimes S_t(-i)\rightarrow R/(I)_{\ge d}
\rightarrow  0
$$ where $\Sym^i(U)=H^0(\mathcal O_{\Lambda}(i))$ is a vector space of homogeneous forms of degree $i$ generated by $U$.
By sheafifying this exact sequence and tensoring $\bigotimes\mathcal O_{\P^{n-t}}(d-1$), we have the surjective morphism
of sheaves
$$
\cdots\longrightarrow \bigoplus^{d-1}_{i=0}\Sym^i(U)\otimes \mathcal O_{\P^{n-t}}(d-1-i)\longrightarrow{\pi_{\Lambda}}_{\ast}\mathcal O_X(d-1)\longrightarrow 0.
$$
For all $y\in Y_t$, we have the following surjective commutative diagram $(\ast)$ by Nakayama's lemma:
\[
\begin{array}{ccccccccccccccccccccccc}
\,\,\, \ds \bigoplus^{d-1}_{i=0}\Sym^i(U)\otimes \mathcal O_{\P^{n-t}}(d-1-i)\otimes k(y)& \ds\rightarrow \,\,\,{\pi_{\Lambda}}_{*}\mathcal O_X(d-1) \otimes k(y) & \ds\rightarrow 0 \\[1ex]
\ds \|&(\ast)\,\,\,\,\,\,\,\,\,\,\,\,\,\,\,\,\,\,\,\,\,\,\,\,\,\,\,\,\,\,\,\,\,\,\,\,\,\,\,\|\, \,\,\,\,\,\,\,\,\,\,\,\,\,\,\,\,\,\,\,\ds & \\[1ex]
\ds H^0(\langle \Lambda, y\rangle, {\mathcal O}_{\langle \Lambda, y\rangle}(d-1))& \ds\rightarrow \,\,\,H^0({\mathcal O}_{{\pi_{\Lambda}}^{-1}(y)}(d-1))& \ds\rightarrow 0 .\\[1ex]
\end{array}
\]
Therefore, as a finite scheme,  ${\pi_{\Lambda}}^{-1}(y)$ is $(d-1)$-normal for all $y\in Y_t$.

For a proof of (b), suppose that $\reg(X\cap L)> d$ for some linear section $X\cap L$ as a finite scheme where $L=\P^{k_0}$
for some  $1\le k_0\le p$. Then we can take a linear subspace $\Lambda_1 \subset L$ of dimension ${k_0}-1$ disjoint from
$X\cap L$. Then $X\cap L$ is a fiber of projection $\pi_{\Lambda_1}:X\to \P^{n-k_0-1}$ at $\pi_{\Lambda_1}(L)$. However,
this is a contradiction by $(a)$.

Let's proceed to prove (c). Suppose that $I$ satisfies $\textbf {N}_{3,p+1}$ and $\textbf {N}_{2,p}$. So $\beta^R_{p+1, 2}$ is nonzero and $\beta^R_{p+1, j}=0$ for
all $j\ge 3$.  Let $S_{p+1}=k[x_{p+1},x_{p+2}\ldots,x_n]\subset S_{p}= k[x_{p}, x_{p+1}\ldots,x_n]$.
Then it follows from Theorem~\ref{Prop:201} that the minimal free presentation of $R/I$ as a $S_{p}$-module is of the form:
\[\cdots \longrightarrow F_1\longrightarrow S_{p}\oplus S_{p}(-1)^{p}\longrightarrow R/I\longrightarrow 0\]
Now consider the following long exact sequence for each $j=0,1,2$,
\begin{equation}\label{eq:210}
\begin{array}{ccccccccccccc}
\st{} \Tor_{1}^{S_p}(R/I,k)_{j}& \st{}
\Tor_{0}^{S_{p+1}}(R/I,k)_{j-1}&\st{\delta= \times x_p}&&
\\[2ex]
&\Tor_{0}^{S_{p+1}}(R/I,k)_{j}&\st{}& \Tor_{0}^{S_p}(R/I,k)_{j}&\st{} &0.
\end{array}
\end{equation}
By the property $\textbf {N}_{3,p+1}$ of $R/I$, we can easily verify that the minimal free resolution of $R/I$ as a $S_{p+1}$-module is of the form
\begin{equation}\label{eq:213}
 \cdots \longrightarrow S_{p+1}\oplus S_{p+1}(-1)^{p+1}\oplus S_{p+1}(-2)^{\alpha}\longrightarrow R/I\longrightarrow 0
\end
{equation}
for some $\alpha$ in $\mathbb Z_{\geq 0}$. Then it follows from ~\ref{eq:210} and $j=2$ that
$$\dim_k\Tor_{0}^{S_{p+1}}(R/I,k)_{1}=p+1\ge \alpha .$$

On the other hand, we have the following surjections from the fact
that
for ${1\le i \le p+1}$, (cf, Proposition ~\ref {Prop:202})
$\Tor_{p+1-i}^{S_{i}}(R/I,k)_{p+1-i+3}=0 $:
\[\Tor_{p+1}^{R}(R/I,k)_{p+1+2}\lra \Tor_{p}^{S_{1}}(R/I,k)_{p+2}\st{\times x_0} \Tor_{p}^{S_{1}}(R/I,k)_{p+3}=0,\]
\[\Tor_{p+1-i}^{S_{i}}(R/I,k)_{p+1-i+2}\lra \Tor_{p-i}^{S_{i+1}}(R/I,k)_{p-i+2}\st{\times x_i}\Tor_{p-i}^{S_{i+1}}(R/I,k)_{p-i+3}=0,\]
\[\Tor_{1}^{S_{p}}(R/I,k)_{3}\lra \Tor_{0}^{S_{p+1}}(R/I,k)_{2}\st{\times x_p} \Tor_{0}^{S_{p+1}}(R/I,k)_{3}= 0\]
for all $0\leq i\leq p+1$, which implies that
\[\dim_k\Tor_{p+1}^{R}(R/I,k)_{p+3}=\beta^R_{p+1,2}(R/I)\ge \alpha=\dim_k\Tor_{0}^{S_{p+1}}(R/I,k)_2.\]
So, we have the inequality
$$\alpha \le \min\{p+1, \beta^R_{p+1,2}(R/I)\}.$$
This completes the proof of (c) by sheafification of the sequence (\ref {eq:213}).
\end{proof}

\begin{Remk} For $p=0$, the bound in $(c)$ is clearly sharp because $X$ is cut out by at most cubic equations. For $p=1$,
we can check that $\alpha=1$. In addition, there is a unique conic passing through general $5$ points,
there is no $5$-secant $2$-plane to $X$ which is cut out by quadrics. In particular,
if $\beta^R_{p+1,2}=0$ then we know that $X$ has no $(p+3)$-secant $(p+1)$-plane because of the property $\textbf{N}_{2, p+1}$.
Therefore, it would also be interesting to check whether the upper bound $2p+3$ in $(c)$ is optimal if $p+1 < \beta^R_{p+1,2}$
for some $p\ge 2$.
\end{Remk}

\begin{Ex}
(a) Let $S^{\ell}(C)$ be the ${\ell}$-th higher secant variety of the rational normal curve $C$ in $\mathbb P^n$. Then the defining ideal
of $S^{\ell}(C)$ is generated by maximal minors of 1-generic matrix of linear forms in $S=\mathbb C[x_0,\ldots,x_n]$ of size ${\ell}+1$.
Then, it follows that $S^{\ell}(C)$ is aCM of degree $\binom{n-{\ell}+1}{{\ell}}$ having $({\ell}+1)$-linear resolution which is given by
Eagon-Northcott complex. Let $\Lambda=\mathbb P^{n-2{\ell}}$ be a general linear space and consider a linear projection from the center $\Lambda$.
Then the length of a general fiber is the degree of $S^{\ell}(C)$ which is equal to the dimension of the space of $\ell$ forms on the linear span of
the fiber. So, the bound in Theorem~\ref{Theo:1-3}~(a) is sharp. For example, if $X$ has 3-linear resolution up to $p$-th step and a linear space $L$ is a $\ell$-secant $p$-plane,
then $\ell\leq \binom{p+2}{2}$, which is sharp as the example $S^2(C)$ in $\mathbb P^5$ shows. \\
(b) Let $C$ be an elliptic normal curve of degree $d$ in $\P^{d-1}$ which satisfies property $\textbf{N}_{d-3}$ but fails to satisfy
property $\textbf{N}_{d-2}$ with $\beta^R_{d-2,2}=1$. Since $\deg(C)=d= d-1 + \min\{d-2, \beta^R_{p+1,2}(R/I)\}$,
the bound in Theorem~\ref{Theo:1-3},(c) is also sharp for the case $\beta^R_{d-2,2}(R/I)< d-2$.\\
(c) Let $C=\nu_{10}(\P^1)$ be a rational normal curve in  $\P^{10}$.
Let $S^{\ell}(C)$ be the $\ell$-th higher secant variety of $\dim S^{\ell}(C)=\min\{2\ell-1, 10\}$.
Then,
$$C\subsetneq \Sec(C)=S^{2}(C)\subsetneq S^{3}(C) \subsetneq \cdots \subsetneq S^{6}(C)=\P^{10}.$$
Then, for any point $q\in S^5(C)\setminus S^{4}(C)$, $\pi_{q}(C)\subset \P^{9}$ is a smooth rational curve with property $\textbf {N}_{2,2}$
(\cite{P1}) with $\beta^R_{3,2}=1$. Since $\pi_{q}(C)$ has a $5$-secant $3$-plane in $\P^{9}$, the bound in Theorem~\ref{Theo:1-3},(c)
is also sharp for this case.
\end{Ex}

\begin{Remk}
In the process of proving Theorem ~\ref{Theo:1-3}, we know that the global property $\textbf {N}_{d,p}$ scheme-theoretically gives local information
on the length of fibers in any linear projection from the center $\Lambda$ of dimension \,$\le p-1$. The commutative diagram $(\ast)$
in the proof can also be understood geometrically as follows:

\[
\begin{array}{cclll}
{\Bl}_{\Lambda}(\mathbb P^n)&= &\mathbb P(\mathcal O_{\mathbb P^{n-t}}(1)\oplus U\otimes \mathcal O_{\mathbb P^{n-t}})&\subset \mathbb P^n \times \mathbb P^{n-t}\\
\xymatrix{
\ar[d]_{\sigma}\\
{}} & &\xymatrix{\ar[d]_{\rho} \\
{}}\\
\mathbb P^n &\xymatrix{ \ar@{.>}[rr]^{\pi_{\Lambda }}&&} &\mathbb P^{n-t}&\\
\end{array}
\]
where $\sigma:\Bl_{\Lambda}{\P^n}\rightarrow \P^n$ is a blow-up of $\P^n$ along $\Lambda$ and
$$\rho:\P(\mathcal O_{\P^{n-t}}(1)\oplus U\otimes\mathcal O_{\P^{n-t}})\rightarrow \P^{n-t}$$ is a vector bundle over $\P^{n-t}$.
We have a natural morphism of sheaves
\[\rho_{\ast}\sigma^{\ast}\mathcal O_{\mathbb P^n}(d-1) \longrightarrow \rho_{\ast}\sigma^{\ast}\mathcal O_X(d-1)= {\pi_{\Lambda}}_{\ast}\mathcal O_X(d-1) \]
where $\mathcal E =\mathcal O_{\P^{n-t}}(1)\oplus U\otimes\mathcal O_{\P^{n-t}}$ and $\rho_{\ast}\sigma^{\ast}\mathcal O_{\mathbb P^n}(d-1)=
\rho_{\ast}\mathcal O_{\mathbb P(\mathcal E)}(d-1)=\Sym^{d-1}(\mathcal E)$. We actually showed from the property $N_{d,p}$ that the following morphism

\[
\begin{array}{rcccccccccccccccccccccc}
\Sym^{d-1}(\mathcal E)\otimes k(y)&\longrightarrow &{\pi_{\Lambda}}_{\ast}\mathcal O_X(d-1)\otimes k(y)\\[2ex]
\parallel\hspace{0.8cm} & & \parallel\\
 \ds \bigoplus^{d-1}_{i=0}\Sym^i(U)\otimes \mathcal O_{\P^{n-t}}(d-1-i)\otimes k(y)& \longrightarrow &H^0({\mathcal O}_{{\pi_{\Lambda}}^{-1}(y)}(d-1))\\[1ex]
\parallel\hspace{0.8cm}& & \\[1ex]
H^0(\langle \Lambda, y\rangle, {\mathcal O}_{\langle \Lambda, y\rangle}(d-1))&&
\end{array}
\]
is surjective for all $y\in Y_t$.
Similar constructions were used in bounding regularity of smooth surfaces and threefolds in \cite{K} and \cite{L}.
\end{Remk}

\section{Effects of property $\textbf{N}_{2,p}$ on projections and multiple loci}\label{section_3}

For a projective variety $X\subset \P^n$, property $\textbf{N}_{2,p}$ is a natural generalization of property $\textbf{N}_p$.
Note that a smooth variety $X\subset \mathbb P^n$ satisfying property $\textbf{N}_{2,p}, p\ge 1$
scheme-theoretically has $\reg(X)\le e+1$ where $e=\codim(X, \P^n)$ and so $X$ is $m$-normal for all $m\ge e$ ~(cf. \cite{BEL}).
The main theorems in this section show that property $\textbf{N}_{2,p}$ plays an important role to control the normality and defining equations
of projected varieties under isomorphic and birational projections  up to $(p-1)$-th step.
\begin{Thm}\label{Main result N_{2,p}}{\bf(Isomorphic projections for $\textbf{N}_{2,p}$ case)}\newline
Suppose that $X\subset \mathbb P^n$ satisfy property
$\textbf{N}_{2,p}$ scheme-theoretically for some $p\geq 2$. Consider an isomorphic projection
$\pi_{\Lambda}:X\to Y_t\subset \mathbb P^{n-t}$ for some $1\leq t\le p-1$. Then we have the following:
\begin{enumerate}
\item[(a)] $H^1(\mathcal I_X(m))=H^1(\mathcal I_{Y_t}(m))$ for all $m\ge t+1$. Thus, if $X$ is
$m$-normal for all $m\geq \alpha_X$ then $Y_t$ is $m$-normal for all $m\ge {\max}\,\,\{\alpha_X, t+1\}$
and $\reg(Y_t)\le \max\{\reg(X), t+2\}$;
\item[(b)] In particular, if $I_X$ satisfies $\textbf{N}_{2,p}$ then $I_{Y_t}$ is also cut out by equations of degree at most $t+2$
and further satisfies property $\textbf {N}_{t+2,p-t}$.
\end{enumerate}
\end{Thm}
\begin{proof}
Let $R=k[x_0,x_1\ldots,x_n]$ and $S_t=k[x_t, x_{t+1},\ldots,x_n]$ be the coordinate rings of $\P^n$ and $\P^{n-t}$ respectively.
Choose an ideal $I$ defining $X$ with $\textbf{N}_{2,p}$ scheme-theoretically. Then, by Theorem ~\ref{Prop:201}, we have the minimal
free resolution of $R/I$ as a graded $S_t$-module:
$$\rightarrow S_t(-p+t-1)^{\oplus \beta_{p-t,1}}\rightarrow \cdots \rightarrow S_t(-2)^{\oplus \beta_{1,1}}
\rightarrow S_t\oplus S_t(-1)^{\oplus \,\,t}\st{\varphi_0} R/I
\rightarrow  0.
$$
Note that $\pi_{{\Lambda}_{*}}(\mathcal O_X)\simeq \mathcal O_{Y_t}$
and $(R/I)_d=H^0(\mathcal O_{Y_t}(d))$ for all $d>>0$.
Therefore, by sheafifying the resolution of $R/I$, we have the following familiar two diagrams by  using Snake Lemma(\cite{GLP},\cite{KP});

\begin{equation}\label{diagram:301}
\begin{array}{ccccccccccccccccccccccc}
&&0&&0&&&&\\[1ex]
&&\downarrow &&\downarrow&&&&\\[1ex]
0&\lra& \mathcal I_{Y_t} & \rightarrow & \mathcal O_{\P^{n-t}}& \rightarrow & \mathcal O_{Y_t}&\rightarrow & 0 \\[1ex]
&&\downarrow && \downarrow && \,\parallel& \\[1ex]
0&\lra & \mathcal L & \rightarrow & {\mathcal O_{\P^{n-t}}(-1)^{\oplus t}}\oplus {\mathcal O_{\P^{n-t}}} & \st{\widetilde{\varphi_0}} &
{\mathcal O_{Y_t}}&\rightarrow & 0\\[1ex]
&&\downarrow &&\downarrow &&  &&\\[1ex]
&&{\mathcal O_{\P^{n-t}}(-1)^{\oplus t}}& = & {\mathcal O_{\P^{n-t}}(-1)^{\oplus t}} &&&& \\[1ex]
&&\downarrow &&\downarrow&&&&\\[1ex]
&&0&&0&&&&
\end{array}
\end{equation}

and in the first syzygies of $R/I$, we have the following diagram:

\begin{equation}\label{diagram:302}
\begin{array}{ccccccccccccccccccccccc}
&&&&&&0&&\\[1ex]
&&&&&&\downarrow&&\\[1ex]
&&0&&&&\mathcal I_{Y_{t}}\\[1ex]
&&\downarrow&&&&\downarrow&&\\[1ex]
0&\lra& \mathcal K & \rightarrow & \mathcal O_{\P^{n-t}}(-2)^{\oplus {\beta}_{1,1}}& \rightarrow &
\mathcal L&\rightarrow & 0 \\[1ex]
&& \downarrow && \,\parallel && \,\downarrow & \\[1ex]
0&\lra & \mathcal N & \rightarrow &  {\mathcal O_{\P^{n-t}}(-2)^{\oplus {\beta}_{1,1}}}& \rightarrow &
{\mathcal O_{\P^{n-t}}(-1)^{\oplus t}}
&\rightarrow & 0\\[1ex]
&&\downarrow &&&&\downarrow  &&\\[1ex]
&&{\mathcal I_{Y_t}}&&&&0&& \\[1ex]
&&\downarrow &&&&&&\\[1ex]
&&0&&&&&&
\end{array}
\end{equation}

\begin{Clm}\label{Claim:101} From the commutative diagrams (\ref{diagram:301}) and (\ref{diagram:302}),
\begin{enumerate}
\item[(a)] $\reg({\mathcal N})\le t+2$;
\item[(b)] For all $m\ge t+1$, we have the isomorphisms on $m$-normality: \\
 $H^1(\mathcal I_{Y_t/\P^{n-t}}(m))\simeq H^2(\mathcal K(m))\simeq H^1(\mathcal L(m))\simeq H^1(\mathcal I_{X/\P^n}(m))$.
 \end{enumerate}
\end{Clm}
\begin{proof}
First of all, we can control the Castelnuovo-regularity of $\mathcal N$ (cf. \cite {GLP}, \cite {KP}, \cite {N})
by using Eagon-Northcott complex associated to the exact sequence
$$0\rightarrow {\mathcal N}\rightarrow {\mathcal O_{\P^{n-t}}(-2)^{\oplus {\beta}_{1,1}}}\rightarrow
{\mathcal O_{\P^{n-t}}(-1)^{\oplus t}}\rightarrow 0.$$ As a consequence,
$\reg({\mathcal N})\le t+2$. Thus, from the leftmost column and first row of (\ref{diagram:302}), we have the following isomorphisms
for all $m\ge t+1$:
$$H^1(\mathcal I_{Y_t/\P^{n-t}}(m))\simeq H^2(\mathcal K(m))\simeq H^1(\mathcal L(m)). $$

On the other hand, by taking global sections and using simple linear syzygies of $R/I$ as $S_t$-module, we have
the following two commutative diagrams:
\[
\begin{array}{ccccccccccccccccccccccc}
0&\lra &\ds H^0_{*}(\mathcal L) & \ds\rightarrow & \ds S_t(-1)^{\oplus \,\,t}\oplus S_t& \ds\st{H^0_*(\widetilde{\varphi_0})} & \ds\oplus_{\ell \in \Z} H^0 (\mathcal{O}_{Y_t}(\ell))& \ds\rightarrow  H^1_{*}(\mathcal L)& \ds\rightarrow 0 \\[1ex]
&& \,\,\,\,\ds\uparrow && \,\,\quad\ds\parallel && \,\,\,\,\,\ds\uparrow &  \\[1ex]
0&\ds\lra &\ds K_0& \ds\rightarrow & \ds S_t(-1)^{\oplus \,\,t}\oplus S_t& \st{\varphi_0} & R/I &\rightarrow 0\,\,.&&
\end{array}
\]
Since $\im\,\, H^0_{*}(\widetilde{\varphi_0})=R/I_X$ and $\oplus_{\ell \in \Z} H^0 (\mathcal{O}_{Y_t}(\ell))=\oplus_{\ell \in \Z} H^0 (\mathcal{O}_{X}(\ell))$,
we get
$H^1_{*}(\mathcal L)=H^1_{*}(\mathcal I_{X/\P^n})$. Therefore, our claim and (a) are proved.\end{proof}

Now, let's return to the proof of (b) in Theorem~\ref{Main result N_{2,p}}. In this case, note that $I=I_X$ and $H^1_{*}(\mathcal K)=0$.
Consider the following diagram for all $\ell\ge 1$:
\[
\begin{array}{ccccccccccccccccccccccc}
\ds H^0({\mathcal O}_{\P^{n-t}}(\ell))\otimes H^0({\mathcal N}(t+2))& \ds\twoheadrightarrow \,\,\,H^0({\mathcal N}(t+2+\ell))& \ds\rightarrow 0 \\[1ex]
\ds\downarrow &\ds\downarrow & \\[1ex]
\ds H^0({\mathcal O}_{\P^{n-t}}(\ell))\otimes H^0(\mathcal I_{Y_t}(t+2))& \ds\twoheadrightarrow \,\,\,H^0({\mathcal I}_{Y_t}(t+2+\ell))& \ds\rightarrow 0 \\[1ex]
\ds\downarrow &\ds\downarrow & \\[1ex]
\ds0 &\ds 0 & \\[1ex]
\end{array}
\]
Note that surjectivity of the first row is given by $\reg(\mathcal N)\le t+2$ and surjectivity of two vertical columns are given
by the fact $H^1_{*}(\mathcal K)=0$. Thus, the second row is also surjective and consequently $Y_t$ is cut out by equations of
degree at most $(t+2)$. For the syzygies of $I_{Y_t}$, consider the exact sequence by taking global sections
$$
0\rightarrow H^0_{*}({\mathcal K})\rightarrow H^0_{*}({\mathcal N})\rightarrow  I_{Y_t}\rightarrow H^1_{*}({\mathcal K})=0.
$$
Since $H^0_{*}(\mathcal K)=K_1$ is the first syzygy module of $R/I_X$, $H^0_{*}(\mathcal K)$ has the following resolution:
$$\rightarrow S_t(-p+t-1)^{\oplus \beta_{p-t,1}}\rightarrow \cdots \rightarrow S_t(-4)^{\oplus \beta_{3,1}}\rightarrow S_t(-3)^{\oplus \beta_{2,1}}\rightarrow H^0_{*}(\mathcal K)\rightarrow  0
$$
and so, $\Tor^{S_t}_i(H^0_{*}(\mathcal K),k)_{i+j}=0$ for $0\le i \le p-t-2$ and $j\ge 4.$
On the other hand, we know the following equivalence;
$$\reg \,H^0_{*}(\mathcal N)=\reg({\mathcal N})\le t+2 \Longleftrightarrow \Tor^{S_t}_i(H^0_{*}(\mathcal N),k)_{i+j}=0 \,\,{\text for}\,\, i\ge 0,j\ge t+3.$$
Thus, from the long exact sequence:
\begin{equation*}
\begin{array}{ccccccccccccc}
\Tor_{i}^{S_t}(H^0_{*}({\mathcal K}),k)_{i+j}\lra\Tor_{i}^{S_t}(H^0_{*}({\mathcal N}),k)_{i+j}\lra\Tor_{i}^{S_t}(I_{Y_t},k)_{i+j}
\lra&  &&
\\[2ex]
\st{\delta}\Tor_{i-1}^{S_t}(H^0_{*}({\mathcal K}),k)_{i+j}\lra\Tor_{i-1}^{S_t}(H^0_{*}({\mathcal N}),k)_{i+j}\lra\Tor_{i-1}^{S_t}(I_{Y_t},k)_{i+j},
\end{array}
\end{equation*}
we get $\Tor^{S_t}_i(I_{Y_t},k)_{i+j}=0$ for $0\le i \le p-t-1$ and $j\ge t+3$ and $Y_t$ satisfies property $N_{2+t, p-t}$.
\end{proof}

In the complete embedding of $X\subset \P(H^0(\mathcal O_X(1)))$, property $\textbf{N}_{2,p}$ is the same as property $\textbf{N}_p$.
In this case, we have the following Corollary which is already given in Theorem 1.2 in \cite {KP} and Corollary 3 in \cite {CKP}.
\begin{Coro}\label{Np-normality}
Let $X\subset \mathbb P(H^0(\mathcal O_X(1)))=\P^n$ be a reduced non-degenerate projective variety with property $\textbf{N}_p$
for some $p\geq 2$. Consider an isomorphic projection $\pi_{\Lambda}:X\to Y_t\subset \P(W)=\P^{n-t}, t=\codim(W, H^0(\mathcal O_X(1))), 1\le t\le p-1$.
The projected variety $Y_t\subset \P(W)$ satisfies the following:
\begin{enumerate}
\item[(a)] $Y_t$ is $m$-normal for all $m\ge t+1$;
\item[(b)] $Y_t$ is cut out by equations of degree at most $(t+2)$ and further satisfies property $\textbf{N}_{2+t,p-t}$;
\item[(c)] $\reg(Y_t)\le \max\{\reg(X), t+2\}$.
\end{enumerate}
\end{Coro}
\begin{proof}
This is clear from Theorem ~\ref{Main result N_{2,p}} with $n_0(X)=1$. For a different proof using vector
bundle technique in the restricted Euler sequence, see \cite {CKP} and \cite {KP}.
\end{proof}

\begin{Remk} A. Alzati and F. Russo  gave a necessary and sufficient condition for the isomorphic projection of a $m$-normal variety
to be also $m$-normal. As an application, they showed that for a variety $X\subset \P^n$ satisfying property $\textbf{N}_2$, one point isomorphic
projection of $X$ in $\P^{n-1}$ is $k$-normal for all $k\ge 2$ (Theorem 3.2 and Corollary 3.3 in [1]). So, Theorem ~\ref{Main result N_{2,p}}
is a generalization to nonlinearly normal embedding on normality, defining equations and their syzygies.
\end{Remk}
\begin{Ex}
Let $C=\nu_{13}(\P^1)$ be a rational normal curve in  $\P^{13}$.
Let $S^{\ell}(C)$ be the $\ell$-th higher secant variety of $\dim S^{\ell}(C)=\min\{2\ell-1,13\}$.
Then, for any point $q\in \P^{13}\setminus S^{6}(C)$, $\pi_{q}(C)\subset \P^{12}$ is a smooth rational curve with property $\textbf{N}_{2,4}$
(\cite{P1}) and for a general line $\ell \subset \P^{13}$, $\pi_{\ell}(C)\subset \P^{11}$ is a rational curve with property $\textbf{N}_{2,3}$
by using Singular or Macaulay 2.
So, by Theorem 3.1.(a), $\pi_{\ell}(C)\subset \P^{11}$ and $\pi_{\Lambda}(C)\subset \P^{10}$ are $m$-normal for all $m\ge 2$
and $3$-regular for general plane $\Lambda$ in $\P^{13}$. This is a refinement of Corollary \ref{Np-normality}.
\end{Ex}


On the other hand, for a point $q\in \Sec(X)\cup \Tan(X)$ we can also consider a birational projection and syzygies of
the projected varieties. To begin with, let us explain the basic situation and information on the partial elimination ideals
under outer projections. For $q=(1,0,\cdots, 0,0)\notin X$, consider an outer projection $\pi_{q}:X\to Y_1\subset \mathbb P^{n-1}=\Proj(S_1), S_1=k[x_1, x_{2},\ldots,x_n]$. Suppose the ideal $I$ define $X$ scheme-theoretically.
For the degree lexicographic order, if $f\in I$ has leading term $\ini(f)=x_0^{d_0}\cdots x_n^{d_n}$, we set $d_0(f)=d_0$,
the leading power of $x_0$ in $f$. Then it is well known that
$K_0(I)=\bigoplus_{m\ge 0}\big\{f \in I_m \mid d_0(f)=0\big\}=I\cap S_1$ and defines $Y_1$ scheme-theoretically.
More generally, let us give the definition and basic properties of partial elimination ideals, which was introduced by M.\ Green in \cite{G}.

\begin{Def}[\cite{G}]\label{def_partial elimi ideals}
Let $I\subset R$ be a homogeneous ideal and let  \[\tilde{K}_i(I)=\bigoplus_{m\ge 0}\big\{f\in I_{m}\mid d_0(f)\leq i\big\}.\]
If $f\in \tilde{K}_i(I)$, we may write uniquely $f=x_0^i\bar{f}+g$ where $d_0(g)<i$. Now we define $K_i(I)$
by the image of $\tilde{K}_i(I)$ in $S_1$ under the map $f\mapsto \bar{f}$ and we call $K_{i}(I)$ the $i$-th
partial elimination ideal of $I$.
\end{Def}

Note that $K_{0}(I)=I\cap S_1$ and there is a short exact sequence as graded $S_1$-modules
\begin{equation}\label{eq:201}
0\rightarrow \frac{\tilde{K}_{i-1}(I)}{\tilde{K}_{0}(I)}
\rightarrow \frac{\tilde{K}_{i}(I)}{\tilde{K}_{0}(I)}\rightarrow
K_{i}(I)(-i)\rightarrow 0.
\end{equation}
In addition, we have the filtration on partial elimination ideals of $I$:
\[K_0(I)\subset K_1(I)\subset K_2(I)\subset \cdots \subset K_i(I)\subset \cdots \subset S_1.\]

\begin{Prop}[\cite{G}]\label{multipleloci}
Set theoretically, the $i$-th partial elimination ideal $K_i(I)$
is the ideal of\,  $Z_{i}=\big\{q \in Y_1\mid \mult_q(Y_1)\ge i+1\big\}$ for every $i\geq 1$.
\end{Prop}

\begin{Lem}\label{Lem bir}
Let $X\subset \mathbb P^n$ be a reduced non-degenerate projective variety satisfying property $\textbf{N}_{2,p}, p\ge 2$ scheme-theoretically.
Consider a projection $\pi_{q}:X\to Y_1\subset \P^{n-1}$ where $q\notin X$. Suppose
$$\Sigma_{q}(X):=\{x\in X|{\pi_q}^{-1}(\pi_{q}(x)) \text { has length}\ge 2\}$$ be the nonempty secant locus of one-point projection. Then,
\begin{itemize}
\item[(a)] $\Sigma_q(X)$ is a quadric hypersurface in a linear subspace $L$ and $q\in L$;
\item[(b)] $\pi_{q}(\Sigma_{q}(X))=Z_1$ is a linear space which is the support of cokernel of $\mathcal O_{Y_1}\hookrightarrow {\pi_q}_{*}(\mathcal O_X)$;
\item[(c)] For a point $q\in \Sec(X)\setminus \Tan(X)\cup X$,  $\Sigma_q(X)$= $\{\text\,{two \,\,distinct \,\,points}\,\}$.
\end{itemize}
\end{Lem}
\begin{proof}
Since $X$ satisfies $\textbf{N}_{2,p}, p\ge 2$, there is no $4$-secant $2$-plane to $X$ by Theorem ~\ref{Theo:1-3}.(b).
Let $Z_1:=\{y\in Y_1|{\pi_q}^{-1}(y) \text { has length}\ge 2\ \}$ and choose two points $y_1, y_2$ in $Z_1$.
Consider the line $\ell =\overline{y_1, y_2}$ in $\P^{n-1}$. If $\langle y_1, y_2 \rangle\cap Y_1$ is finite, then we have
$4$-secant plane $\langle q, y_1, y_2\rangle$ which is a contradiction. So, $\Sec(Z_1)=Z_1$ and finally, we conclude that $Z_1$
is a linear space. Since $\pi_{q}:\Sigma_q(X)\twoheadrightarrow Z_1\subset Y_1$ is a 2:1 morphism, $\Sigma_q(X)$ is a quadric hypersurface
in $L=\langle Z_1, q\rangle$.  For a proof of (c), if $\dim \Sigma_q(X)$ is positive, then clearly, $q\in \Tan\,\Sigma_q(X)\subset \Tan(X)$.
So, we are done. \end{proof}

We have also the following generalization of the Lemma ~\ref{Lem bir} for $\textbf{N}_{d,2}, d\ge 2$ by using the mapping cone theorem
and partial elimination ideals..
\begin{Prop}\label{N_{d,p}-geometry}
Let $X\subset \mathbb P^n$ be a non-degenerate reduced projective scheme and the ideal $I$ define $X$ scheme-theoretically with property $\textbf{N}_{d,2}$.
For any projection $\pi_q: X \to Y=\pi_q(X)\subset \mathbb P^{n-1}$ from a point $q\in \P^n \setminus X$,
$$Z_{i}=\{ y\in \pi_q(X)\,|\, \pi_q^{-1}(y)\text{ has length }\geq i+1\,\}$$
satisfies the following properties for $d-2\leq i\leq d$.
\begin{itemize}
\item [(a)] $K_{d-1}(I)$ is generated by at most linear forms. Thus $Z_{d-1}$ is either empty or a linear space;
\item [(b)] $K_{d-2}(I)$ is generated by at most cubic forms. Thus $Z_{d-2}$ is either empty or cut out by at most cubic equations set-theoretically.
\end{itemize}
\end{Prop}
\begin{proof}
(a): Since the ideal $I$ satisfies property $\textbf{N}_{d,2}$, there exists the following exact sequence~(not necessarily minimal) by Proposition~\ref{Prop:201}:
$$\rightarrow\cdots \rightarrow \bigoplus^{d-1}_{j=1}S_1(-1-j)^{\beta_{1,j}^{S_1}} \st{\varphi_1} \bigoplus^{d-1}_{i=0}S_1(-i) \st{\varphi_0}R/I \rightarrow  0.
$$
Furthermore, we can easily verify that
$\ker \varphi_0$ is $\tilde{K}_{d-1}(I)$ and thus we have the following exact sequence:
\begin{equation}\label{eq:401}
0\rightarrow \tilde{K}_{d-1}(I) \rightarrow \oplus_{i=0}^{d-1}
S_1(-i) \rightarrow R/I\rightarrow 0.
\end{equation}

Now consider the following commutative diagram with  $K_{0}(I)=I\cap S_1$:
\begin{equation}\label{diagram:302-1}
\begin{array}{ccccccccccccccccccccccc}
&&0&&0&&0&&\\[1ex]
&&\downarrow &&\downarrow&&\downarrow&&\\[1ex]
0&\rightarrow& K_{0}(I) & \rightarrow & S_1& \rightarrow &
S_1/K_{0}(I)&\rightarrow & 0 \\[1ex]
&&\downarrow && \downarrow && \,\,\,\downarrow{\tilde \alpha} & \\[1ex]
0&\rightarrow& \tilde{K}_{d-1}(I) & \rightarrow &
\oplus_{i=0}^{d-1} S_1(-i)& \st{\varphi_0} & R/I
&\rightarrow & 0\\[1ex]
&&\downarrow &&\downarrow && \downarrow &&\\[1ex]
0&\rightarrow& {\tilde{K}_{d-1}(I)}/K_{0}(I) & \rightarrow &
\oplus_{i=1}^{d-1} S_1(-i)& \rightarrow& \coker\,\, \tilde \alpha &\rightarrow &0\\[1ex]
&&\downarrow &&\downarrow&&\downarrow&&\\[1ex]
&&0&&0&&0&&
\end{array}
\end{equation}

Since $I$ satisfies the property $\textbf{N}_{d,2}$, it follows from the middle row and left column sequences in the diagram~(\ref{diagram:302-1})
that $\tilde{K}_{d-1}(I)$ and $\tilde{K}_{d-1}(I)/{K_0(I)}$ are generated by at most degree $d$ elements.

On the other hands, we have a short exact sequence from (\ref{eq:201}) :
\begin{equation}\label{diagram:partial}
0\rightarrow \frac{\tilde{K}_{d-2}(I)}{K_0(I)} \rightarrow \frac{\tilde{K}_{d-1}(I)}{K_0(I)}\rightarrow
K_{d-1}(I)(-d+1)\rightarrow 0,
\end{equation}
Hence, $K_{d-1}(I)$ is generated by at most linear forms and further $1$-regular. So, $Z_{d-1}$ is either empty or a linear space
by Proposition~\ref{multipleloci}.

(b): Since $K_{d-1}(I)$ is $1$-regular, we have
$$\Tor_1^{S_1}(K_{d-1}(I)(-d+1), k)_j=0 \,\,\,\text {for all} \,\,\,j\ge d+2$$ and it follows from (\ref{diagram:partial}) that
${\tilde{K}_{d-2}(I)}/{K_0(I)}$ is generated by at most degree $d+1$ elements.
Similarly, consider again the following short exact sequence as $S_1$-modules
\begin{equation*}
0\rightarrow \frac{\tilde{K}_{d-3}(I)}{K_0(I)}\rightarrow \frac{\tilde{K}_{d-2}(I)}{K_0(I)}\rightarrow
K_{d-2}(I)(-d+2)\rightarrow 0.
\end{equation*}
Hence, $K_{d-2}(I_X)$ is generated by at most cubic forms and therefore we complete the proof of (b).
\end{proof}

\begin{Coro}\label{center and multiple locus}
In the same situation as in Proposition ~\ref{N_{d,p}-geometry}, assume that $X$ satisfies property $\textbf{N}_{d, 2}, d\ge 2$.
If the linear space $Z_{d-1}(X)$ is nonempty, then $\langle Z_{d-1}(X), q\rangle\cap X$ is a hypersurface of degree $d$
in the span $\langle Z_{d-1}(X), q\rangle$.
\end{Coro}
\begin{proof}
It is cleat that $\langle Z_{d-1}(X), q\rangle\cap X$ is a hypersurface in a linear space $\langle Z_{d-1}(X), q\rangle$.
Since there is no $d+1$-secant line through $q$, we are done.
\end{proof}

As shown in Lemma ~\ref{Lem bir}, the fact that $Z_1$ is a linear space is crucial in the proof of the following theorem.
\begin{Thm}\label{birational N_{2,p}}{\bf(Birational projections for $\textbf{N}_{2,p}$ case)}\newline
Let $X\subset \mathbb P^n$ be a reduced non-degenerate projective variety satisfying property
$\textbf{N}_{2,p}$ scheme-theoretically for $p\geq 2$. Consider a projection
$\pi_{q}:X\to Y_1\subset\mathbb P^{n-1}$ where $q\in \Sec(X)\cup \Tan(X)\setminus X$.  Then we have the following:
\begin{itemize}
\item[(a)] $H^1_*(\mathcal I_X)=H^1_*(\mathcal I_{Y_1})$. Thus, $Y_1$ is $m$-normal if and only if $X$ is $m$-normal
for all $m\ge 1$, and $\reg(Y_1)\le \max\{\reg(X), \reg(\mathcal O_{Y_1})+1\}$;
\item[(b)] $Y_1$ is cut out by at most cubic hypersurfaces and satisfies property $N_{3,p-1}$.
\end{itemize}
\end{Thm}
\begin{proof} We may assume that $q=(1,0,,,0)\in \Sec(X)\cup \Tan(X)\setminus X$.
Let $R=k[x_0,x_1\ldots,x_n]$ be a coordinate ring of $\P^n$, $S_1=k[x_1, x_{2},\ldots,x_n]$ be a coordinate ring of $\P^{n-1}$.
Let the ideal $I$ define $X$ with the condition $\textbf{N}_{2,p}$ scheme-theoretically. Then,
it is easily checked that $K_0(I)$ also defines $Y_1$ scheme-theoretically.
By Theorem ~\ref{Prop:201}, we have the minimal free resolution of $R/I$ as a graded $S_1$-module:
$$\cdots\rightarrow S_1(-p)^{\oplus \beta_{p-1,1}}\rightarrow\cdots \rightarrow S_1(-2)^{\oplus \beta_{1,1}}
\st{\varphi_1} S_1\oplus S_1(-1) \st{\varphi_0} R/I\rightarrow  0.
$$
Furthermore, we have the following diagram:
\begin{equation*}
\begin{array}{ccccccccccccccccccccccc}
&&0&&0&&0&&\\[1ex]
&&\downarrow &&\downarrow&&\downarrow&&\\[1ex]
0&\rightarrow& K_0(I) & \rightarrow & S_1& \rightarrow &
S_1/K_0(I)&\rightarrow & 0 \\[1ex]
&&\downarrow && \downarrow && \,\,\,\downarrow{\alpha} & \\[1ex]
0&\rightarrow& \tilde{K}_{1}(I) & \rightarrow & S_1\oplus S_1(-1)& \st{\varphi_0} & R/I
&\rightarrow & 0\\[1ex]
&&\downarrow &&\downarrow && \downarrow &&\\[1ex]
0&\rightarrow& K_1(I)(-1) & \rightarrow & S_1(-1)&\rightarrow& \coker\,\, \alpha &\rightarrow &0\\[1ex]
&&\downarrow &&\downarrow&&\downarrow&&\\[1ex]
&&0&&0&&0&&
\end{array}
\end{equation*}

Note that $\varphi_0(f,g)=f+g\cdot x_0$ and thus, $K_1(I)$ is the first partial elimination ideal of $I$ associated
to the projection $\pi_{q}$.
Since $\tilde{K}_{1}(I)$ has the following minimal free resolution as a graded $S_1$-module:
$$\cdots\rightarrow S_1(-p)^{\oplus \beta_{p-1,1}}\st{\varphi_{p-1}}\cdots \rightarrow S_1(-2)^{\oplus \beta_{1,1}}
\st{\varphi_1} \tilde{K}_{1}(I)\rightarrow  0,
$$
we know that $K_1(I)$ is generated by linear forms and
$$\reg(K_1(I)(-1))=2,  \,\,\,\,\coker\,\,\alpha = S_1/K_1(I)(-1).$$
Moreover, by usual $\Tor$-computations, $K_0(I)$ satisfies property $N_{3, p-1}$.

On the other hand, consider the following exact sequence
\begin{equation}\label{coker1}
0\rightarrow \mathcal O_{Y_1}\st{\tilde \alpha}{\pi_q}_{*}(\mathcal O_X)\rightarrow\coker\,\, \tilde \alpha\rightarrow 0.
\end{equation}
By Lemma~\ref{Lem bir}, since $\coker \,\,\tilde \alpha$ has the support $Z_1$ which is a linear space in $\P^{n-1}$ and
 $\pi_{q}:\Sigma_q(X)\twoheadrightarrow Z_1$ is $2:1$, we have
 $${\pi_q}_{*}(\mathcal O_X)|_{Z_1}=\mathcal O_{Z_1}\oplus\mathcal O_{Z_1}(-1) \,\text { and }\,\,\, \coker\,\,\tilde \alpha= \mathcal O_{Z_1}(-1).$$
Therefore, $H^0_{*}(\coker\,\,\alpha)=S_1/I_{Z_1}(-1)$.
Then, by taking global sections from the above sequence ~(\ref{coker1}), we have the following commutative diagram as $S_1$-modules
with exact rows and columns:
\begin{equation*}
\begin{array}{ccccccccccccccccccccccc}
&&0&&0&&0&&\\[1ex]
&&\downarrow &&\downarrow&&\downarrow&&\\[1ex]
0&\rightarrow& S_1/I_{Y_1} & \rightarrow & \bigoplus_{m=0}^{\infty}H^0(\mathcal O_{Y_1}(m))& \rightarrow &
H^1_{*}(\mathcal I_{Y_1})&\rightarrow & 0 \\[1ex]
&&\downarrow && \downarrow && \,\downarrow& \\[1ex]
0&\lra & R/I_X & \rightarrow & \bigoplus_{m=0}^{\infty}H^0(\mathcal O_X(m)) & \rightarrow & H^1_{*}(\mathcal I_X)
&\rightarrow & 0\\[1ex]
&&\downarrow &&\downarrow && \downarrow &&\\[1ex]
&& S_1/I_{Z_1}(-1)& =  & S_1/I_{Z_1}(-1)&&0&&\\[1ex]
&&\downarrow &&\downarrow&&&&\\[1ex]
&&0&&0&&&&
\end{array}
\end{equation*}
The reason why the left column is exact is that $K_1(I)=I_{Z_1}$ for any ideal $I$ defining $X$ scheme-theoretically.
Thus, $H^1_{*}(\mathcal I_{Y_1})\simeq H^1_{*}(\mathcal I_X)$ and so, $X$ is $m$-normal if and only if $Y_1$ is $m$-normal.
So we complete the proof of (a) and (b).
\end{proof}

\begin{Remk}
For a complete embedding of $X\subset \P(H^0(\mathcal O(1))$ satisfying property $N_p$, Lemma ~\ref{Lem bir} and
Theorem~\ref{birational N_{2,p}} was proved in \cite{P2} with different method. However, the point is that we can also deal
with non-complete embeddings of $X$ in $\P^n$ satisfying property $\textbf{N}_{2,p}$ by virtue of the graded mapping cone theorem
without using Green-Lazarsfeld's vector bundle technique on restricted Euler sequence on $X$.
\end{Remk}

\section{ Moving the center and the structure of projected varieties}\label{section_4}

In the previous sections, we proved the uniform properties of higher normality and syzygies of projections
when the given variety $X$ satisfies property $\textbf {N}_{2,p}, p\ge 2$ scheme-theoretically.
However, according to moving the center, we have a lot of interesting varieties with different structures in algebra,
geometry and syzygies. As an example, for a rational normal curve $C=\nu_d(\P^1)$ in $\P^d$, consider
the following filtration on the $\ell$-th higher secant variety $S^{\ell}(C)$ of dimension $\min\{2{\ell}-1, d\}$:
$$
C\subsetneq \Sec(C)=S^{2}(C)\subsetneq S^{3}(C) \subsetneq \cdots \subsetneq S^{\lfloor\frac{d}{2}\rfloor}(C)\subsetneq
S^{\lfloor \frac{d}{2}\rfloor+1}(C)=\P^d\,\,\,,$$
Then we have (\cite{CKP}, \cite{P1})
\begin{enumerate}
\item[(a)] $\overline{\pi_{q}(C)}\subset \P^{d-1}$ satisfies property $\textbf {N}_{2, d-2}$ for $q\in C$,
\item[(b)] $\pi_{q}(C)\subset \P^{d-1}$ is a rational curve with one node satisfying property $\textbf {N}_{2, d-3}$ for
$q\in \Sec(C)\setminus C$,
\item[(c)] $\pi_{q}(C)\subset \P^{d-1}$ satisfies property $\textbf {N}_{2, \ell-3}$ for $q\in S^{\ell}(C)\setminus S^{\ell-1}(C)$.
\end{enumerate}
Note that all projected curves are $m$-normal for all $m\ge 2$ and thus $3$-regular.

Note that for varieties of next to minimal degree, the arithmetic properties of projected varieties by moving the center
were investigated in \cite{BS} for the first time.
The following proposition show that the number of quadratic equations, Hilbert functions, the depth of projected varieties
and Betti tables depend on the dimension of the secant locus $\Sigma_{q}(X)$ and the position of the center of projection.
For a complete embedding $X\subset \P(H^0(\mathcal L))$,
the same result is given in \cite{P2}. Let $s=\dim \Sigma_{q}(X)$ and if the secant locus $\Sigma_{q}(X)=\emptyset$, then $s=-1$.

\begin{Prop}\label{prop_moving}
Let $X\subset \P^n$ be a reduced non-degenerate projective variety satisfying property $\textbf {N}_{2,p},\, p\ge 2$. Consider the projection
$\pi_{q}:X\to Y_1\subset\mathbb P^{n-1}$ where $q\notin X$.
Let $\Sigma_{q}(X)$ is the secant locus of the projection $\pi_{q}$. Then the following holds:
\begin{enumerate}
\item[(a)]
$h^0(\P^{n-1}, \mathcal I_{Y_1}(2))= h^0(\P^{n}, \mathcal I_{X}(2))-n+s$,
\item[(b)]
$\depth(Y_1)=\min\{\depth(X), s+2\}$  under the condition that
$$H^i(\mathcal O_X(j))=0, \forall j\le -i, 1\le i \le \dim(X).$$
\end{enumerate}
\end{Prop}
\begin{proof}
(a) First, for the isomorphic projection case, we obtained the following fact from the commutative diagram (3.2):
$$
\reg(\mathcal N)=3,\,\,\,\, h^1(\mathcal I_{Y_1}(\ell))=h^2(\mathcal K(\ell))=h^1(\mathcal L(\ell))=h^1(\mathcal I_X(\ell))\,\,\,\text{for}\,\,\, \ell\ge 2.
$$ From the basic equalities
\begin{equation*}
                  \left\{
                  \begin{array}{ll}
                  h^0(\mathcal I_{X}(2))+ h^0(\mathcal O_X(2))=\binom {n+2}{2}+ h^1(\mathcal I_X(2)) \,\,\,\text {and}\,\,\,\\[1ex]
                  h^0(\mathcal I_{Y_1}(2))+ h^0(\mathcal O_{Y_1}(2))=\binom{n+1}{2}+ h^1(\mathcal I_{Y_1}(2)),
                  \end{array}\right.\\
\end{equation*}
we get $h^0(\mathcal I_{Y_1}(2))= h^0(\mathcal I_{X}(2))-n-1$. In the case of finite birational projections,
the secant locus $\Sigma_{q}(X)$ is not empty and $\pi_{q}(\Sigma_{q}(X))=Z_1=\P^s$.
In the proof of Theorem ~\ref{birational N_{2,p}}, we got the following fact:
$$H^1_{*}(\P^{n-1}, \mathcal I_{Y_1})\simeq H^1_{*}(\P^n, \mathcal I_X),\,\,\, 0\rightarrow S_1/I_{Y_1}\rightarrow R/I_X\rightarrow S_1/I_{Z_1}(-1)\rightarrow 0.
$$
Therefore, by simple computation we have $h^0(\mathcal I_{Y_1}(2))= h^0(\mathcal I_{X}(2))-n+s$.
For a proof of (b), consider the following exact sequence
\begin{equation}\label{coker2}
0\rightarrow \mathcal O_{Y_1}\st{\alpha}{\pi_q}_{*}(\mathcal O_X)\rightarrow\mathcal O_{Z_1}(-1)\rightarrow 0.
\end{equation}
If $s=-1$, then $Z_1=\emptyset$ and $H^1(\mathcal I_{Y_1}(1))\neq 0$. So, $\depth(Y_1)=1$.
Suppose $\depth(X)=1, s\ge 0$. Then, by Theorem~\ref{birational N_{2,p}} (a), $H^1_{*}(\mathcal I_{Y_1})\simeq H^1_{*}(\mathcal I_X)\neq 0$ and $\depth(Y_1)=1$.

Now, suppose $\depth(X)\ge 2, s\ge 0$. When $s=0$, then $Z_1$ is one point. Therefore we have
$$0\rightarrow H^0(\mathcal O_{Z_1}(\ell-1))\rightarrow H^1(\mathcal O_{Y_1}(\ell))\rightarrow H^1(\mathcal O_{X}(\ell))\rightarrow 0$$ and $0\neq H^0(\mathcal O_{Z_1}(\ell-1))\subset H^1(\mathcal O_{Y_1}(\ell))$ for all $\ell\le 0$.
So, $\depth(Y_1)=2=\min\{\depth(X), s+2\}$.  For $s\ge 1$ and $\depth(X)\ge s+2$, we obtain the sequence
$$
0\rightarrow H^s_{*}(\mathcal O_{Z_1}(-1))\rightarrow H^{s+1}_{*}(\mathcal O_{Y_1})\rightarrow H^{s+1}_{*}(\mathcal O_{X})\rightarrow 0,
$$
$H^{i}_{*}(\mathcal O_{Y_1})\simeq H^{i}_{*}(\mathcal O_{X})=0$ for all $1\le i\le s$ and  $H^s_{*}(\mathcal O_{Z_1}(-1))\neq 0$.
Thus, $\depth(Y_1)=s+2=\min\{\depth(X), s+2\}$. Finally, in the case of $2\le\depth(X)\le s+1, s\ge 1,$ under the assumption that $H^i(\mathcal O_X(j))=0, \forall j\le -i$, we can easily check that $\depth(Y_1)=\depth(X)=\min\{\depth(X), s+2\}$. Therefore, we are done.
\end{proof}
\begin{Remk}
Similary, under the same assumption, we have
$$h^0(\P^{n-1}, \mathcal I_{Y_1}(d))= h^0(\P^{n}, \mathcal I_{X}(d))-\binom{n+d-1}{n}+\binom{s+d-1}{s}$$
which shows that Hilbert functions of projected varieties depend only on the dimension of the singular locus.
But, the Betti tables are much more delicate and we get the following additional information:
for each $1\leq i\leq p-1$,
\begin{equation*}\label{eq:406}
\beta_{i,2}(S/I_{Y_1})-\beta_{i+1,1}(S/I_{Y_1})=\binom{n-s-1}{i+1}+(-1)^i+
\sum_{1\leq j\leq i+1}(-1)^{j+i}\beta_{j,1}^R(R/I_X).
\end{equation*}
In particular, for $i=1$,
\begin{equation*}\label{eq:407} \beta_{1,2}(S/I_{Y_1})-\beta_{2,1}(S/I_{Y_1})=
\binom{n-s-1}{2}-1+\beta_{1,1}^R(R/I_X)- \beta_{2,1}^R(R/I_X).
\end{equation*}
Hence, we see that if $\binom{n-s-1}{2}-1+\beta_{1,1}^R(R/I_X) - \beta_{2,1}^R(R/I_X)$ is positive then there exist cubic generators of $I_{Y_1}$.
\end{Remk}
\begin{Ex}(A non-normal variety with non-vanishing cohomology)\\
We give some examples related to our proposition.
For a projective normal variety $X$, we define $$\delta(X):=\min\{\depth \,\mathcal O_{X, x}|\, x \,\,\text {is a closed point}\}.$$
Then $H^i(\mathcal O_X(\ell))=0$ for all $\ell << 0$ and $i < \delta(X)$ by vanishing theorem of Enriques-Severi-Zariski-Serre.
In the proof of proposition ~\ref{prop_moving}, for $s=0$ we have an interesting example $Y_1$ such that $Y_1$ has only one isolated
non-normal singular point and in fact, $H^1(\mathcal O_{Y_1}(\ell))\neq 0$ for all $\ell \le 0$.
As examples, suppose that a projective variety $X$ has no lines and plane conics in $\P^n$ with the condition $\textbf {N}_{2,p}, p\ge 2$
(e.g., the Veronese variety $\upsilon_d(\P^n), d\ge 3$ or its isomorphic projections). Then, the singular locus of any simple projection
is either empty or only one point because the secant locus is a quadric hypersurface in some linear subspace.
\end{Ex}

%
%

\begin{Ex}
Consider a rational normal 3-fold scroll $S_{1,1,4}=\P(\mathcal O_{\P^1}(1)\oplus \mathcal O_{\P^1}(1)\oplus\mathcal O_{\P^1}(4))$ in $\P^8$.
From the Eagon-Northcott complex, we obtain the minimal free resolution of $S_{1,1,4}$ as follows:
\[0\rightarrow R(-6)^{5}\rightarrow  R(-5)^{24}\rightarrow R(-4)^{45} \rightarrow  R(-3)^{40} \rightarrow  R(-2)^{15}\rightarrow I_{S_{1,1,4}} \rightarrow 0.\]
As the center of projection $q\in \P^8$ moves toward $S_{1,1,4}$, we will see that
the number of cubic generators decreases and the number of quadric generators increases in the following:
\begin{enumerate}
\item[(a)]
Let $q\in \P^8\setminus \Sec(S_{1,1,4})$ and any isomorphic projection $Y \subset \P^7$ has
the following resolution with $\depth(Y)=1$
$$
\cdots \rightarrow  S(-4)^{40}\oplus S(-3)^{8} \rightarrow  S(-3)^{10}\oplus S(-2)^{6}\rightarrow I_Y \rightarrow 0.
$$
\item[(b)]
Suppose $q\in \Sec(S_{1,1,4})\setminus \Tan(S_{1,1,4})$.
Then $s=0$ and $I_Y$ has the following resolution with $\depth(Y)=2$:
$$
\dots \rightarrow  S(-4)^{19}\oplus S(-3)^{8} \rightarrow  S(-3)^{3}\oplus S(-2)^{7}\rightarrow I_Y \rightarrow 0.
$$
\item[(c)] For a point $q\in \Tan(S_{1,1,4})\setminus S_{1,1,4}$, $Y$ has different two types of resolutions:
First, consider a linear span $\P^3 = \langle \ell_1, F\rangle$ where $\ell_1$ is a line embedded by
$\P(\mathcal O_{\P^1}(1))\hookrightarrow\P(\mathcal O_{\P^1}(1)\oplus \mathcal O_{\P^1}(1)\oplus\mathcal O_{\P^1}(4))\subset \P^8$
and $F$ be the any fiber of the morphism $\varphi:S_{1,1,4}\rightarrow \P^1$. For a general point $q\in \P^3= \langle \ell_1, F\rangle$,
$Y$ has a singular locus $\P^1$, only one cubic generator and the following minimal resolution of length $5$:
$$
\cdots \rightarrow S(-4)^4\oplus S(-3)^{12} \rightarrow  S(-3)\oplus S(-2)^{8} \rightarrow I_Y \rightarrow 0.$$
Second, take a general point $q\in \P^3$ where the quadric hypersurface $\P(\mathcal O_{\P^1}(1)\oplus \mathcal O_{\P^1}(1))\subset \P^3$
is a subvariety of $S_{1,1,4}\subset \P^8$. Then the projected variety $Y$ has clearly the singular locus $\P^2$, $\depth(Y)=4$ and
the following resolution:
\[0\rightarrow S(-6)\rightarrow S(-4)^{9} \rightarrow  S(-3)^{16} \rightarrow  S(-2)^{9} \rightarrow I_Y \rightarrow 0.\]
\item[(d)]
For a general point $q\in S_{1,1,4}$, an inner projection $Y$ is a smooth 3-fols scroll of type $S_{1,1,3}$ and  has the following resolution:
\[0\rightarrow S(-5)^4\rightarrow S(-4)^{15} \rightarrow  S(-3)^{20} \rightarrow  S(-2)^{10} \rightarrow I_Y \rightarrow 0.\]
\end{enumerate}
\end{Ex}

As mentioned in \cite {EGHP1}, property $\textbf{N}_{2,p}$ is {\it rigid} : if $X$ is a reduced subscheme in $\P^n$ with the condition
$\textbf{N}_{2,p}, p=\codim(X, \P^n)$, then $X$ is $2$-regular.
By the rigidity, for a projective reduced scheme $X$ satisfying property $\textbf {N}_{2,p}, p=\codim(X, \P^n)$, any outer projection $\pi_{q}(X)$ does not satisfy
the property $\textbf {N}_{2,p-1}$.
On the other hand, for a projective variety $X\subset \P^n$ with the condition $\textbf {N}_{2,p}, p\ge 2$
and $q\notin X$, we obtained that $\pi_{q}(X)$ satisfies at least property $\textbf {N}_{3, p-1}$ by Theorem~\ref{Main result N_{2,p}} and Theorem~\ref{birational N_{2,p}}.

However, we raise the following question for inner projections:
\begin{Qu}
Let $X$ be a projective reduced scheme in $\P^n$ satisfying property $\textbf {N}_{2,p}, p\ge 1$ which is not necessarily linearly normal.
Consider the inner projection from linear subvariety $L$ of $X$ and $Y=\overline{\pi_{L}(X\setminus L)}$ in $\P^{n-t-1}$, where $\dim L=t < p$.
In contrast with the outer projections, is it true that $Y$ satisfies $\textbf {N}_{2, p-t-1}$ for a linear space $L\subset X$?
For example, for a nondegenerate smooth variety $X$ in $\P(H^0(\mathcal L))$ with property $\textbf {N}_p$, $Y$ satisfies $\textbf {N}_{p-1}$ for a point
$q\in X\setminus \Trisec(X)$ where $\Trisec(X)$ is the union of all proper trisecant lines or lines in $X$~(see \cite{CKK} for details).
Note that the graded mapping cone theorem can not be directly applied to this case because $R/I_X$ is infinitely generated as a $S_t$-module.
\end{Qu}

\bibliographystyle{amsalpha}

\end{document}